\documentclass[11pt,twoside]{amsart}
\usepackage{graphicx}
\usepackage{amssymb}
\usepackage{amsmath}
\usepackage{amsthm}
\usepackage{amsfonts}
\usepackage{amsmath, amsfonts, amssymb}
\newtheorem{theorem}{Theorem}[section]
\newtheorem{corollary}[theorem]{Corollary}

\newtheorem{definition}[theorem]{Definition}
\newtheorem{remark}[theorem]{\bf{Remark}}

\numberwithin{equation}{section}

\begin{document}
\begin{center}\small{In the name of Allah, Most Gracious, Most Merciful.}\end{center}
\vspace{1cm}

\title{Classifications of two-dimensional Jordan algebras over the algebraically closed fields and $\mathbb{R}$}
\author{H.Ahmed$^1$, U.Bekbaev$^2$, I.Rakhimov$^3$}

\thanks{{\scriptsize
emails: $^1$houida\_m7@yahoo.com; $^2$bekbaev@iium.edu.my; $^3$rakhimov@upm.edu.my.}}
\maketitle
\begin{center}
\address{$^{1}$Department of Math., Faculty of Science, UPM, Selangor, Malaysia $\&$ Depart. of Math., Faculty of Science, Taiz University, Taiz, Yemen}\\
  \address{$^2$Department of Science in Engineering, Faculty of Engineering, IIUM, Malaysia}\\
\address{$^3$Institute for Mathematical Research (INSPEM), UPM, Serdang, Selangor, Malaysia}
\end{center}

\begin{abstract} The classification, up to isomorphism, of two-dimensional (not necessarily commutative) Jordan algebras over algebraically closed fields and $\mathbb{R}$ is presented in terms of their matrices of structure constants.
\end{abstract}
{\scriptsize Keywords: basis; Jordan algebra; matrix of structural constants.

}
\pagestyle{myheadings}
\markboth{\rightline {\sl   H. Ahmed, U.Bekbaev, I.Rakhimov}}
         {\leftline{\sl   Classification of two-dimensional real Jordan algebras}}

\bigskip
\section{Introduction}

The classification problem of finite dimensional algebras is one of the important problems of algebra. Many authors have considered this problem for different classes of two-dimensional algebras \cite{A,A1, Bermúdez, J, Bu, C, D, M, M1} and in \cite{P2000} the problem is considered for all two-dimensional algebras over any field. The approach of \cite{P2000} is basis free (invariant) approach. Unlike the approach of \cite{P2000} a complete list of isomorphism classes of two-dimensional algebras over algebraically closed fields is presented in \cite{B} by the use of structure constants approach and the same approach to the problem over the  field of real numbers $\mathbb{R}$ has been applied in \cite{Bekbaev}. 

One of the interesting classes of algebras is the class of Jordan algebras.
In \cite{Bermúdez, Wallace}\ and \cite{M} the authors determined the isomorphic classes of Jordan algebras (commutative case) over the field of real numbers and over algebraically closed fields, respectively.

In this paper we make use results of \cite{B} and \cite{Bekbaev} to give complete lists of two-dimensional (not necessarily commutative) Jordan algebras over any algebraically closed and real fields by providing lists of canonical representatives of their matrices of structure constants. To the best of the authors' knowledge such  classifications of two-dimensional Jordan algebras has not been given yet. For some computations we used Computer Software Mathematica.

\section{Preliminaries}

Let $\mathbb{F}$ be any field and $A\otimes B$ stand for the Kronecker product that is for the matrix with blocks $(a_{ij}B),$ where $A=(a_{ij})$ and $B$ are matrices over $\mathbb{F}$.

\begin{definition} A vector space  $\mathbb{A}$ over $\mathbb{F}$ with multiplication $\cdot :\mathbb{A}\otimes \mathbb{A}\rightarrow \mathbb{A}$ given by $(\mathbf{u},\mathbf{v})\mapsto \mathbf{u}\cdot \mathbf{v}$ such that \[(\alpha\mathbf{u}+\beta\mathbf{v})\cdot \mathbf{w}=\alpha(\mathbf{u}\cdot \mathbf{w})+\beta(\mathbf{v}\cdot \mathbf{w}),\ \ \mathbf{w}\cdot (\alpha\mathbf{u}+\beta\mathbf{v})=\alpha(\mathbf{w}\cdot \mathbf{u})+\beta(\mathbf{w}\cdot \mathbf{v})\] whenever $\mathbf{u}, \mathbf{v}, \mathbf{w}\in \mathbb{A}$ and $\alpha, \beta\in \mathbb{F}$, is said to be an algebra.\end{definition}

\begin{definition} Two algebras $\mathbb{A}$ and $\mathbb{B}$ are called isomorphic if there is an invertible linear map  $\mathbf{f}:\mathbb{A}\rightarrow \mathbb{B} $ such that \[\mathbf{f}(\mathbf{u}\cdot_{\mathbb{A}} \mathbf{v})=\mathbf{f}(\mathbf{u})\cdot_{\mathbb{B}} \mathbf{f}(\mathbf{v})\] whenever $\mathbf{u}, \mathbf{v}\in \mathbb{A}.$\end{definition}

Let ($\mathbb{A}, \cdot)$ be a $m$-dimensional algebra over $\mathbb{F}$ and $e=(e^1,e^2,...,e^m)$ be its basis. Then the
bilinear map $\cdot$ is represented by a matrix $A=(A^k_{ij})\in M(m\times m^2;\mathbb{F})$ as follows \[\mathbf{u}\cdot
\mathbf{v}=eA(u\otimes v),\] for $\mathbf{u}=eu,\mathbf{v}=ev,$
where $u = (u_1, u_2, ..., u_m)^T,$ and $v = (v_1, v_2, ..., v_m)^T$ are column coordinate vectors of $\mathbf{u}$ and
$\mathbf{v},$ respectively.
The matrix $A\in M(m\times m^2;\mathbb{F})$ defined above is called the matrix of structural constants (MSC) of
$\mathbb{A}$ with respect to the basis $e$. Further we assume that a basis $e$ is fixed and we do not make a difference
between the algebra
$\mathbb{A}$ and its MSC $A$.

If $e'=(e'^1,e'^2,...,e'^m)$ is another basis of $\mathbb{A}$, $e'g=e$ with $g\in G=GL(m;\mathbb{F})$, and  $A'$ is MSC
of $\mathbb{A}$ with respect to $e'$ then it is known that
\[A'=gA(g^{-1})^{\otimes 2}\] is valid (see \cite{B}). Thus, we can reformulate the isomorphism of algebras as follows.

\begin{definition} Two $m$-dimensional algebras $\mathbb{A}$, $\mathbb{B}$ over $\mathbb{F}$, given by
	their matrices of structure constants $A$, $B$, are said to be isomorphic if $B=gA(g^{-1})^{\otimes 2}$ holds true
	for some $g\in GL(m;\mathbb{F})$.\end{definition}

Further we consider only the case $m=2$ and for the simplicity we use $A=\left(\begin{array}{cccc} \alpha_1 & \alpha_2 &
\alpha_3 &\alpha_4\\ \beta_1 & \beta_2 & \beta_3 &\beta_4\end{array}\right)$ for MSC, where $\alpha_1, \alpha_2,
\alpha_3, \alpha_4, \beta_1, \beta_2, \beta_3, \beta_4$ stand for any elements of $\mathbb{F}$.

Due to \cite{B} we have the following classification theorems according to $Char(\mathbb{F})\neq 2,3,$ $Char(\mathbb{F})=2$ and $Char(\mathbb{F})=3$ cases, respectively.
\begin{theorem}\label{T1} Over an algebraically closed field $\mathbb{F}$ $(Char(\mathbb{F})\neq 2$ and $3)$, any non-trivial $2$-dimensional algebra is isomorphic to only one of the following algebras listed by their matrices of structure constants:
	\begin{itemize}
		\item $A_{1}(\mathbf{c})=\left(
		\begin{array}{cccc}
		\alpha_1 & \alpha_2 &\alpha_2+1 & \alpha_4 \\
		\beta_1 & -\alpha_1 & -\alpha_1+1 & -\alpha_2
		\end{array}\right),\ \mbox{where}\ \mathbf{c}=(\alpha_1, \alpha_2, \alpha_4, \beta_1)\in \mathbb{F}^4,$
		\item $A_{2}(\mathbf{c})=\left(
		\begin{array}{cccc}
		\alpha_1 & 0 & 0 & 1 \\
		\beta _1& \beta _2& 1-\alpha_1&0
		\end{array}\right)\simeq \left(
		\begin{array}{cccc}
		\alpha_1 & 0 & 0 & 1 \\
		-\beta _1& \beta _2& 1-\alpha_1&0
		\end{array}\right),\ \mbox{where}\ \mathbf{c}=(\alpha_1, \beta_1, \beta_2)\in \mathbb{F}^3,$
		\item $A_{3}(\mathbf{c})=\left(
		\begin{array}{cccc}
		0 & 1 & 1 & 0 \\
		\beta _1& \beta _2 & 1&-1
		\end{array}\right),\ \mbox{where}\ \mathbf{c}=(\beta_1, \beta_2)\in \mathbb{F}^2,$
		\item $A_{4}(\mathbf{c})=\left(
		\begin{array}{cccc}
		\alpha _1 & 0 & 0 & 0 \\
		0 & \beta _2& 1-\alpha _1&0
		\end{array}\right),\ \mbox{where}\ \mathbf{c}=(\alpha_1, \beta_2)\in \mathbb{F}^2,$
		\item $A_{5}(\mathbf{c})=\left(
		\begin{array}{cccc}
		\alpha_1& 0 & 0 & 0 \\
		1 & 2\alpha_1-1 & 1-\alpha_1&0
		\end{array}\right),\ \mbox{where}\ \mathbf{c}=\alpha_1\in \mathbb{F},$
		\item $A_{6}(\mathbf{c})=\left(
		\begin{array}{cccc}
		\alpha_1 & 0 & 0 & 1 \\
		\beta _1& 1-\alpha_1 & -\alpha_1&0
		\end{array}\right)\simeq \left(
		\begin{array}{cccc}
		\alpha_1 & 0 & 0 & 1 \\
		-\beta _1& 1-\alpha_1 & -\alpha_1&0
		\end{array}\right),\ \mbox{where}\ \mathbf{c}=(\alpha_1, \beta_1)\in \mathbb{F}^2,$
		\item $A_{7}(\mathbf{c})=\left(
		\begin{array}{cccc}
		0 & 1 & 1 & 0 \\
		\beta_1& 1& 0&-1
		\end{array}\right),\ \mbox{where}\ \mathbf{c}=\beta_1\in \mathbb{F},$
		\item $A_{8}(\mathbf{c})=\left(
		\begin{array}{cccc}
		\alpha_1 & 0 & 0 & 0 \\
		0 & 1-\alpha_1 & -\alpha_1&0
		\end{array}\right),\ \mbox{where}\ \mathbf{c}=\alpha_1\in \mathbb{F},$
		\item $A_{9}=\left(
		\begin{array}{cccc}
		\frac{1}{3}& 0 & 0 & 0 \\
		1 & \frac{2}{3} & -\frac{1}{3}&0
		\end{array}\right),$
		\item $A_{10}=\left(
		\begin{array}{cccc}
		0 & 1 & 1 & 0 \\
		0 &0& 0 &-1
		\end{array}
		\right),$
		\item $A_{11}=\left(
		\begin{array}{cccc}
		0 & 1 & 1 & 0 \\
		1 &0& 0 &-1
		\end{array}
		\right),$
		\item $A_{12}=\left(
		\begin{array}{cccc}
		0 & 0 & 0 & 0 \\
		1 &0&0 &0\end{array}
		\right).$
\end{itemize}\end{theorem}

\begin{theorem}\label{T2} Over an algebraically closed field $\mathbb{F}$ $(Char(\mathbb{F})=2)$, any non-trivial $2$-dimensional algebra is isomorphic to only one of the following algebras listed by their matrices of structure constants:
	\begin{itemize}
		\item $A_{1,2}(\mathbf{c})=\left(
		\begin{array}{cccc}
		\alpha_1 & \alpha_2 &\alpha_2+1 & \alpha_4 \\
		\beta_1 & -\alpha_1 & -\alpha_1+1 & -\alpha_2
		\end{array}\right),\ \mbox{where}\ \mathbf{c}=(\alpha_1, \alpha_2, \alpha_4, \beta_1)\in \mathbb{F}^4,$
		\item $A_{2,2}(\mathbf{c})=\left(
		\begin{array}{cccc}
		\alpha_1 & 0 & 0 & 1 \\
		\beta _1& \beta_2 & 1-\alpha_1&0
		\end{array}\right),\ \mbox{where}\ \mathbf{c}=(\alpha_1, \beta_1, \beta_2)\in \mathbb{F}^3,$
		\item $A_{3,2}(\mathbf{c})=\left(
		\begin{array}{cccc}
		\alpha_1 & 1 & 1 & 0 \\
		0& \beta_2 & 1-\alpha_1&1
		\end{array}\right),\ \mbox{where}\ \mathbf{c}=(\alpha_1, \beta_2)\in \mathbb{F}^2,$
		\item $A_{4,2}(\mathbf{c})=\left(
		\begin{array}{cccc}
		\alpha _1 & 0 & 0 & 0 \\
		0 & \beta_2 & 1-\alpha _1&0
		\end{array}\right),\ \mbox{where}\ \mathbf{c}=(\alpha_1,\beta_2)\in \mathbb{F}^2,$
		\item $A_{5,2}(\mathbf{c})=\left(
		\begin{array}{cccc}
		\alpha_1 & 0 & 0 & 0 \\
		1 & 1 & 1-\alpha_1&0
		\end{array}\right),\ \mbox{where}\ \mathbf{c}=\alpha_1\in \mathbb{F},$
		\item $A_{6,2}(\mathbf{c})=\left(
		\begin{array}{cccc}
		\alpha_1 & 0 & 0 & 1 \\
		\beta _1& 1-\alpha_1 & -\alpha_1&0
		\end{array}\right),\ \mbox{where}\ \mathbf{c}=(\alpha_1, \beta_1)\in \mathbb{F}^2,$
		\item $A_{7,2}(\mathbf{c})=\left(
		\begin{array}{cccc}
		\alpha_1 & 1 & 1 & 0 \\
		0& 1-\alpha_1& -\alpha_1&-1
		\end{array}\right),\ \mbox{where}\ \mathbf{c}=\alpha_1\in \mathbb{F},$
		\item $A_{8,2}(\mathbf{c})=\left(
		\begin{array}{cccc}
		\alpha_1 & 0 & 0 & 0 \\
		0 & 1-\alpha_1 & -\alpha_1&0
		\end{array}\right),\ \mbox{where}\ \mathbf{c}=\alpha_1\in \mathbb{F},$
		\item $A_{9,2}=\left(
		\begin{array}{cccc}
		1 & 0 & 0 & 0 \\
		1 & 0 & 1&0
		\end{array}\right),$
		\item $A_{10,2}=\left(
		\begin{array}{cccc}
		0 & 1 & 1 & 0 \\
		0 &0& 0 &-1
		\end{array}
		\right),$
		\item $A_{11,2}=\left(
		\begin{array}{cccc}
		1 & 1 & 1 & 0 \\
		0 &-1& -1 &-1
		\end{array}
		\right),$
		\item $A_{12,2}=\left(
		\begin{array}{cccc}
		0 & 0 & 0 & 0 \\
		1 &0&0 &0\end{array}
		\right).$
	\end{itemize}
\end{theorem}

\begin{theorem}\label{T3} Over an algebraically closed field $\mathbb{F}$ $(Char(\mathbb{F})=3)$, any non-trivial $2$-dimensional algebra is isomorphic to only one of the following algebras listed by their matrices of structure constant matrices:
	\begin{itemize}
		\item $A_{1,3}(\mathbf{c})=\left(
		\begin{array}{cccc}
		\alpha_1 & \alpha_2 &\alpha_2+1 & \alpha_4 \\
		\beta_1 & -\alpha_1 & -\alpha_1+1 & -\alpha_2
		\end{array}\right),\ \mbox{where}\ \mathbf{c}=(\alpha_1, \alpha_2, \alpha_4, \beta_1)\in \mathbb{F}^4,$
		\item $A_{2,3}(\mathbf{c})=\left(
		\begin{array}{cccc}
		\alpha_1 & 0 & 0 & 1 \\
		\beta _1& \beta _2& 1-\alpha_1&0
		\end{array}\right)\simeq\left(
		\begin{array}{cccc}
		\alpha_1 & 0 & 0 & 1 \\
		-\beta _1& \beta _2& 1-\alpha_1&0
		\end{array}\right),\ \mbox{where}\ \mathbf{c}=(\alpha_1, \beta_1, \beta_2)\in \mathbb{F}^3,$
		\item $A_{3,3}(\mathbf{c})=\left(
		\begin{array}{cccc}
		0 & 1 & 1 & 0 \\
		\beta _1& \beta _2 & 1&-1
		\end{array}\right),\ \mbox{where}\ \mathbf{c}=(\beta_1, \beta_2)\in \mathbb{F}^2,$
		\item $A_{4,3}(\mathbf{c})=\left(
		\begin{array}{cccc}
		\alpha _1 & 0 & 0 & 0 \\
		0 & \beta _2& 1-\alpha _1&0
		\end{array}\right),\ \mbox{where}\ \mathbf{c}=(\alpha_1, \beta_2)\in \mathbb{F}^2,$
		\item $A_{5,3}(\mathbf{c})=\left(
		\begin{array}{cccc}
		\alpha_1& 0 & 0 & 0 \\
		1 & -1-\alpha_1 & 1-\alpha_1&0
		\end{array}\right),\ \mbox{where}\ \mathbf{c}=\alpha_1\in \mathbb{F},$
		\item $A_{6,3}(\mathbf{c})=\left(
		\begin{array}{cccc}
		\alpha_1 & 0 & 0 & 1 \\
		\beta _1& 1-\alpha_1 & -\alpha_1&0
		\end{array}\right)\simeq \left(
		\begin{array}{cccc}
		\alpha_1 & 0 & 0 & 1 \\
		-\beta _1& 1-\alpha_1 & -\alpha_1&0
		\end{array}\right),\ \mbox{where}\ \mathbf{c}=(\alpha_1, \beta_1)\in \mathbb{F}^2,$
		\item $A_{7,3}(\mathbf{c})=\left(
		\begin{array}{cccc}
		0 & 1 & 1 & 0 \\
		\beta_1& 1& 0&-1
		\end{array}\right),\ \mbox{where}\ \mathbf{c}=\beta_1\in \mathbb{F},$
		\item $A_{8,3}(\mathbf{c})=\left(
		\begin{array}{cccc}
		\alpha_1 & 0 & 0 & 0 \\
		0 & 1-\alpha_1 & -\alpha_1&0
		\end{array}\right),\ \mbox{where}\ \mathbf{c}=\alpha_1\in \mathbb{F},$
		\item $A_{9,3}=\left(
		\begin{array}{cccc}
		0 & 1& 1& 0 \\
		1 &0&0 &-1\end{array}
		\right),$
		\item $A_{10,3}=\left(
		\begin{array}{cccc}
		0 & 1 & 1 & 0 \\
		0 &0&0 &-1\end{array}
		\right),$
		\item $A_{11,3}=\left(
		\begin{array}{cccc}
		1 & 0 & 0 & 0 \\
		1 &-1&-1 &0\end{array}
		\right),$
		\item $A_{12,3}=\left(
		\begin{array}{cccc}
		0 &0 &0 & 0 \\
		1 &0& 0 &0
		\end{array}
		\right).$
\end{itemize}\end{theorem}
\begin{remark} In \cite{B} the class $A_{3,2}(\mathbf{c})$ should be understood as it is in this paper, as far as there is a type-mistake in this case in \cite{B}.\end{remark}
\section{\bf Classification of $2$-dimensional Jordan algebras over algebraically closed fields}

Recall that an algebra $\mathbb{A}$ is said to be a Jordan algebra if $(\mathbf{u}\cdot
\mathbf{v})\cdot\mathbf{u}^2=\mathbf{u}
\cdot(\mathbf{v}\cdot \mathbf{u}^2)$ for all $\mathbf{u}, \mathbf{v}\in \mathbb{A}$. In terms of its MSC $A$ the
condition above is written as follows:
\[(A(A\otimes A)-A(E\otimes A(E\otimes A)))(u\otimes e_i\otimes u\otimes u)=0,\ \ i=1,2,\]
for $u=\left(
\begin{array}{c}
x \\
y
\end{array}
\right),\ \ e_1=\left(
\begin{array}{c}
1 \\
0
\end{array}
\right),\ \ e_2=\left(
\begin{array}{c}
0 \\
1
\end{array}
\right),\ \ E=\left(
\begin{array}{cc}
1 & 0 \\
0 & 1
\end{array}
\right).$

The condition above in terms of $x,y$ and $\alpha_1, \alpha_2, \alpha_3, \alpha_4, \beta_1, \beta_2, \beta_3, \beta_4$ is
given as a system of equations as follows.

$\begin{array}{c}
x^3 \alpha _1 \alpha _2 \beta _1+x^2 y \alpha _2^2 \beta _1-x^3 \alpha _1 \alpha _3 \beta _1-x^2 y \alpha _3^2 \beta
_1+x^2 y \alpha _1 \alpha _4 \beta _1+2 x y^2 \alpha _2 \alpha _4 \beta _1+y^3 \alpha _4^2 \beta _1 -x^3 \alpha _4 \beta
_1^2+ \\x^3 \alpha _2 \beta _1 \beta _2+x^2 y \alpha _2 \beta _2^2 +x y^2 \alpha _4 \beta _2^2-x^2 y \alpha _1
\alpha _3 \beta _3-x y^2 \alpha _2 \alpha _3 \beta _3-x y^2 \alpha _3^2 \beta _3-y^3 \alpha _3 \alpha _4 \beta _3- \\ 2 x^2 y
\alpha _4 \beta _1 \beta _3+x^2 y \alpha _2 \beta _2 \beta _3-x y^2 \alpha _4 \beta _3^2-x y^2 \alpha _4 \beta _1 \beta
_4+x y^2 \alpha _2 \beta _2 \beta _4+y^3 \alpha _4 \beta _2 \beta _4 -y^3 \alpha _4 \beta _3 \beta _4 =0,\end{array}$\\\\

$\begin{array}{c}-x^2 y \alpha _1 \alpha _3 \beta _1-x y^2 \alpha _2 \alpha _3 \beta _1-x y^2 \alpha _3^2 \beta _1-y^3 \alpha _3 \alpha _4
\beta _1+x^3 \alpha _2 \beta _1^2+2 x^2 y \alpha _2 \beta _1 \beta _2+x y^2 \alpha _4 \beta _1 \beta _2-\\ x^2 y \alpha _1
\beta _2^2-x y^2 \alpha _3 \beta _2^2+x^3 \beta _1 \beta _2^2+x^2 y \beta _2^3 +x^2 y \alpha _1^2 \beta _3+x
y^2 \alpha _1 \alpha _2 \beta _3+x y^2 \alpha _1 \alpha _3 \beta _3+\\ y^3 \alpha _1 \alpha _4 \beta _3-x^3 \alpha _1 \beta
_1 \beta _3+ x^2 y \alpha _2 \beta _1 \beta _3-x^2 y \alpha _3 \beta _1 \beta _3-x y^2 \alpha _4 \beta _1 \beta _3-x^2 y
\alpha _1 \beta _2 \beta _3+x y^2 \alpha _2 \beta _2 \beta _3-\\ x y^2 \alpha _3 \beta _2 \beta _3+x^2 y \beta_2^2 \beta _3-x^2 y \alpha _1 \beta _3^2- 
x y^2 \alpha _3 \beta _3^2-y^3 \alpha _4 \beta _3^2+x^2 y \alpha _1 \beta_1
\beta _4+\\2 x y^2 \alpha _2 \beta _1 \beta _4+x y^2 \alpha _3 \beta _1 \beta _4+y^3 \alpha _4 \beta _1 \beta _4-x^3 \beta _1^2 \beta _4-x y^2 \alpha _1 \beta _2 \beta _4- \\ y^3 \alpha _3 \beta _2 \beta _4+2 x y^2 \beta _2^2 \beta
_4+y^3 \alpha _2 \beta _3 \beta _4-2 x^2 y \beta _1 \beta _3 \beta _4-x y^2 \beta _3^2 \beta _4-x y^2 \beta _1 \beta_4^2+y^3 \beta _2 \beta _4^2-y^3 \beta _3 \beta _4^2=0,\end{array}$\\\\

$\begin{array}{c}-x^3 \alpha _1^2 \alpha _2-x^2 y \alpha _1 \alpha _2^2+x^3 \alpha _1^2 \alpha _3+2 x^2 y \alpha _1 \alpha _3^2+x y^2
\alpha _2 \alpha _3^2+x y^2 \alpha _3^3-x^2 y \alpha _1^2 \alpha _4-2 x y^2 \alpha _1 \alpha _2 \alpha _4+\\ y^3 \alpha _3^2
\alpha _4-y^3 \alpha _1 \alpha _4^2-x^3 \alpha _2^2 \beta _1+x^3 \alpha _1 \alpha _4 \beta _1-x^2 y \alpha _2
\alpha _4 \beta _1+x^2 y \alpha _3 \alpha _4 \beta _1-x^2 y \alpha _2^2 \beta _2-\\ x^3 \alpha _1 \alpha _3 \beta _2-x^2 y
\alpha _2 \alpha _3 \beta _2-x^2 y \alpha _3^2 \beta _2+x^2 y \alpha _1 \alpha _4 \beta _2-x y^2 \alpha _2 \alpha _4
\beta _2-x^3 \alpha _4 \beta _1 \beta _2-x^2 y \alpha _4 \beta _2^2+\\ x^3 \alpha _1 \alpha _2 \beta _3+x^2 y
\alpha _2 \alpha _3 \beta _3+2 x^2 y \alpha _1 \alpha _4 \beta _3+x y^2 \alpha _2 \alpha _4 \beta _3+2 x y^2 \alpha _3
\alpha _4 \beta _3+y^3 \alpha _4^2 \beta _3-\\ x^2 y \alpha _4 \beta _2 \beta _3-x y^2 \alpha _2^2 \beta _4-x^2 y \alpha _1
\alpha _3 \beta _4-x y^2 \alpha _2 \alpha _3 \beta _4-x y^2 \alpha _3^2 \beta _4+x y^2 \alpha _1 \alpha _4
\beta _4-\\ y^3 \alpha _2 \alpha _4 \beta _4+x^3 \alpha _2 \beta _1 \beta _4+x^2 y \alpha _2 \beta _2 \beta _4-x y^2 \alpha
_4 \beta _2 \beta _4+x^2 y \alpha _2 \beta _3 \beta _4+x y^2 \alpha _2 \beta _4^2=0,\end{array}$\\\\

$\begin{array}{c}-x^3 \alpha _1 \alpha _2 \beta _1-x^2 y \alpha _2^2 \beta _1+x^3 \alpha _1 \alpha _3 \beta _1+x^2 y \alpha _3^2 \beta
_1-x^2 y \alpha _1 \alpha _4 \beta _1-2 x y^2 \alpha _2 \alpha _4 \beta _1-y^3 \alpha _4^2 \beta _1+x^3 \alpha _4 \beta
_1^2- \\ x^3 \alpha _2 \beta _1 \beta _2-x^2 y \alpha _2 \beta _2^2-x y^2 \alpha _4 \beta _2^2+x^2 y \alpha _1
\alpha _3 \beta _3+x y^2 \alpha _2 \alpha _3 \beta _3+x y^2 \alpha _3^2 \beta _3+y^3 \alpha _3 \alpha _4 \beta _3+ \\  2 x^2 y
\alpha _4 \beta _1 \beta _3-x^2 y \alpha _2 \beta _2 \beta _3+x y^2 \alpha _4 \beta _3^2+x y^2 \alpha _4 \beta _1 \beta
_4-x y^2 \alpha _2 \beta _2 \beta _4-y^3 \alpha _4 \beta _2 \beta _4+y^3 \alpha _4 \beta _3 \beta _4=0.
\end{array}$\\\\

Equating the coefficients at $x^3, x^2y, xy^2, y^3$ to zero we get the system of equations with respect to\\ $\alpha_1, \alpha_2, \alpha_3, \alpha_4, \beta_1, \beta_2, \beta_3, \beta_4$:
\begin{equation}\label{SE}\begin{array}{rrrrr}
\alpha _1 \alpha _2 \beta _1-\alpha _1 \alpha _3 \beta _1-\alpha _4 \beta _1^2+\alpha _2 \beta _1 \beta _2&=0 \\
\alpha _2^2 \beta _1-\alpha _3^2 \beta _1+\alpha _1 \alpha _4 \beta _1+\alpha _2 \beta _2^2-\alpha _1 \alpha _3 \beta
_3-2 \alpha _4 \beta _1 \beta _3+
\alpha _2 \beta _2 \beta _3 &=0\\
2 \alpha _2 \alpha _4 \beta _1+\alpha _4 \beta _2^2-\alpha _2 \alpha _3 \beta _3-\alpha _3^2 \beta _3-\alpha _4 \beta
_3^2-\alpha _4 \beta _1 \beta _4+
\alpha _2 \beta _2 \beta _4 &=0\\
\alpha _4^2 \beta _1-\alpha _3 \alpha _4 \beta _3+\alpha _4 \beta _2 \beta _4-\alpha _4 \beta _3 \beta _4&=0 \\
\alpha _2 \beta _1^2+\beta _1 \beta _2^2-\alpha _1 \beta _1 \beta _3-\beta _1^2 \beta _4&=0 \\
-\alpha _1 \alpha _3 \beta _1+2 \alpha _2 \beta _1 \beta _2-\alpha _1 \beta _2^2+\beta _2^3+\alpha _1^2 \beta _3+\alpha
_2 \beta _1 \beta _3-\alpha _3 \beta _1 \beta _3 -\alpha _1 \beta _2 \beta _3+\\ \beta _2^2 \beta _3- \alpha _1 \beta
_3^2+\alpha _1 \beta _1 \beta _4-2 \beta _1 \beta _3 \beta _4&=0 \\
-\alpha _2 \alpha _3 \beta _1-\alpha _3^2 \beta _1+\alpha _4 \beta _1 \beta _2-\alpha _3 \beta _2^2+\alpha _1 \alpha _2
\beta _3+\alpha _1 \alpha _3 \beta _3-\alpha _4 \beta _1 \beta _3+\alpha _2 \beta _2 \beta _3-\\ \alpha _3 \beta _2
\beta _3-\alpha _3 \beta _3^2+2 \alpha _2 \beta _1 \beta _4+\alpha _3 \beta _1 \beta _4- \alpha _1 \beta _2 \beta _4+2
\beta _2^2 \beta _4-\beta _3^2 \beta _4-\beta _1 \beta _4^2 &=0\\
-\alpha _3 \alpha _4 \beta _1+\alpha _1 \alpha _4 \beta _3-\alpha _4 \beta _3^2+\alpha _4 \beta _1 \beta _4-\alpha _3
\beta _2 \beta _4+\alpha _2 \beta _3 \beta _4+ \beta _2 \beta _4^2-\beta _3 \beta _4^2&=0 \\
-\alpha _1^2 \alpha _2+\alpha _1^2 \alpha _3-\alpha _2^2 \beta _1+\alpha _1 \alpha _4 \beta _1-\alpha _1 \alpha _3 \beta
_2-\alpha _4 \beta _1 \beta _2+\alpha _1 \alpha _2 \beta _3+\alpha _2 \beta _1 \beta _4&=0 \\
-\alpha _1 \alpha _2^2+2 \alpha _1 \alpha _3^2-\alpha _1^2 \alpha _4-\alpha _2 \alpha _4 \beta _1+\alpha _3 \alpha _4
\beta _1-\alpha _2^2 \beta _2- \alpha _2 \alpha _3 \beta _2- \alpha _3^2 \beta _2+\\ \alpha _1 \alpha _4 \beta _2-
\alpha _4 \beta _2^2+\alpha _2 \alpha _3 \beta _3+2 \alpha _1 \alpha _4 \beta _3-\alpha _4 \beta _2 \beta _3-\alpha _1
\alpha _3 \beta _4+\alpha _2 \beta _2 \beta _4+\alpha _2 \beta _3 \beta _4&=0 \\
\alpha _2 \alpha _3^2+\alpha _3^3-2 \alpha _1 \alpha _2 \alpha _4-\alpha _2 \alpha _4 \beta _2+\alpha _2 \alpha _4 \beta
_3+2 \alpha _3 \alpha _4 \beta _3-\alpha _2^2 \beta _4-\alpha _2 \alpha _3 \beta _4-\\ \alpha _3^2 \beta _4+ \alpha _1
\alpha _4 \beta _4-\alpha _4 \beta _2 \beta _4+\alpha _2 \beta _4^2&=0 \\
\alpha _3^2 \alpha _4-\alpha _1 \alpha _4^2+\alpha _4^2 \beta _3-\alpha _2 \alpha _4 \beta _4&=0.
\end{array}
\end{equation}

\begin{theorem}{\label{thm1}}
	Over an algebraically closed field $\mathbb{F}$, ($Char(\mathbb{F})\neq 2,\ 3$ and $5$), any nontrivial $2$-dimensional
	Jordan algebra is isomorphic to only one of the following algebras listed by their matrices of structure constants, where $i\in\mathbb{F}$
	stands for a fixed element for which $i^2=-1$ :
	\begin{itemize}
		\item $J_1=A_2\left(\frac{1}{2},0,\frac{1}{2}\right)=\left(
		\begin{array}{cccc}
		\frac{1}{2} & 0 & 0 & 1 \\
		0 & \frac{1}{2} & \frac{1}{2} & 0
		\end{array}
		\right),$
		\item $J_2=A_2\left(\frac{1}{2},0,-\frac{1}{2}\right)=\left(
		\begin{array}{cccc}
		\frac{1}{2} & 0 & 0 & 1 \\
		0 & -\frac{1}{2} & \frac{1}{2} & 0
		\end{array}
		\right),$
		\item $J_3(\alpha _1)=A_4(\alpha_1,-1+2\alpha_1)=\left(
		\begin{array}{cccc}
		\alpha _1 & 0 & 0 & 0 \\
		0 & -1+2 \alpha _1 & -\alpha _1+1 & 0
		\end{array}
		\right),$\\  \hspace*{\fill} where $\alpha _1\in \mathbb{F}\ \mbox{and}\ \alpha _1\neq\frac{1}{10} \left(5\pm\sqrt{5}\right),$
				\item $J_4(\alpha
		_1)=A_4\left(\alpha_1,\sqrt{\alpha _1-\alpha _1^2}\right)=\left(
		\begin{array}{cccc}
		\alpha _1 & 0 & 0 & 0 \\
		0 & \sqrt{\alpha _1-\alpha _1^2} & -\alpha _1+1 & 0
		\end{array}
		\right),\ \mbox{where } \alpha _1\in \mathbb{F},$
		\item $J_5(\alpha _1)=A_4\left(\alpha_1,-\sqrt{\alpha _1-\alpha _1^2}\right)=\left(
		\begin{array}{cccc}
		\alpha _1 & 0 & 0 & 0 \\
		0 & -\sqrt{\alpha _1-\alpha _1^2} & -\alpha _1+1 & 0
		\end{array}
		\right),$\\ \hspace*{\fill} where  $\alpha _1\in \mathbb{F}$ and $\alpha _1\neq 0,1,$\\
		\item $J_6=A_5\left(\frac{1}{10}
		\left(5-\sqrt{5}\right)\right)=\left(
		\begin{array}{cccc}
		\frac{1}{10} \left(5-\sqrt{5}\right) & 0 & 0 & 0 \\
		1 & -\frac{\sqrt{5}}{5} & \frac{1}{10} \left(5+\sqrt{5}\right) & 0
		\end{array}
		\right),$
		\item $J_7=A_5\left(\frac{1}{10} \left(5+\sqrt{5}\right)\right)=\left(
		\begin{array}{cccc}
		\frac{1}{10} \left(5+\sqrt{5}\right) & 0 & 0 & 0 \\
		1 & \frac{\sqrt{5}}{5} & \frac{1}{10} \left(5-\sqrt{5}\right) & 0
		\end{array}
		\right),$
		\item $J_8=A_8\left(\frac{1}{3}\right)=\left(
		\begin{array}{cccc}
		\frac{1}{3} & 0 & 0 & 0 \\
		0 & \frac{2}{3} & -\frac{1}{3} & 0
		\end{array}
		\right),$
		\item $J_9=A_8\left(\frac{1}{2}-\frac{i}{2}\right)=\left(
		\begin{array}{cccc}
		\frac{1}{2}-\frac{i}{2} & 0 & 0 & 0 \\
		0 & \frac{1}{2}+\frac{i}{2} & -\frac{1}{2}+\frac{i}{2} & 0
		\end{array}
		\right),$
		\item$J_{10}=A_8\left(\frac{1}{2}+\frac{i}{2}\right)=\left(
		\begin{array}{cccc}
		\frac{1}{2}+\frac{i}{2} & 0 & 0 & 0 \\
		0 & \frac{1}{2}-\frac{i}{2} & -\frac{1}{2}-\frac{i}{2} & 0
		\end{array}
		\right),$
		\item $J_{11}=A_{12}=\left(
		\begin{array}{cccc}
		0 & 0 & 0 & 0 \\
		1 & 0 & 0 & 0
		\end{array}
		\right).$
	\end{itemize}
\end{theorem}
\textbf{Proof}. To prove we make use the above presented results of \cite{B} and test the system of equations (\ref{SE}) for a
consistency.

For $A_{1}(\alpha_1, \alpha_2, \alpha_4, \beta_1)=\left(
\begin{array}{cccc}
\alpha _1 & \alpha _2 & \alpha _2+1 & \alpha _4 \\
\beta _1 & -\alpha _1 & -\alpha _1+1 & -\alpha _2
\end{array}
\right)$ the system of equations (\ref{SE}) becomes\\
\begin{equation}\label{SE1}
\begin{array}{cccc}
-\alpha _1 \beta _1-\alpha _1 \alpha _2 \beta _1-\alpha _4 \beta _1^2 &=0\\
-\alpha _1+\alpha _1^2-2 \alpha _1 \alpha _2+3 \alpha _1^2 \alpha _2-\beta _1-2 \alpha _2 \beta _1-2 \alpha _4 \beta
_1+3 \alpha _1 \alpha _4 \beta _1 &=0\\
-1+\alpha _1-3 \alpha _2+3 \alpha _1 \alpha _2-2 \alpha _2^2+3 \alpha _1 \alpha _2^2-\alpha _4+2 \alpha _1 \alpha _4+3
\alpha _2 \alpha _4 \beta _1 &=0\\
-\alpha _4+\alpha _1 \alpha _4+\alpha _1 \alpha _2 \alpha _4+\alpha _4^2 \beta _1 &=0\\
-\alpha _1 \beta _1+2 \alpha _1^2 \beta _1+2 \alpha _2 \beta _1^2 &=0\\
-\alpha _1+5 \alpha _1^2-6 \alpha _1^3-\beta _1+2 \alpha _2 \beta _1-6 \alpha _1 \alpha _2 \beta _1 &=0\\
-1+4 \alpha _1-4 \alpha _1^2+2 \alpha _1 \alpha _2-6 \alpha _1^2 \alpha _2-\beta _1-4 \alpha _2 \beta _1-6 \alpha _2^2
\beta _1-\alpha _4 \beta _1 &=0\\
-\alpha _1 \alpha _2-2 \alpha _2^2-\alpha _4+3 \alpha _1 \alpha _4-2 \alpha _1^2 \alpha _4-\alpha _4 \beta _1-2 \alpha
_2 \alpha _4 \beta _1 &=0\\
2 \alpha _1^2+\alpha _1 \alpha _2-2 \alpha _2^2 \beta _1+2 \alpha _1 \alpha _4 \beta _1 &=0\\
3 \alpha _1+\alpha _2+7 \alpha _1 \alpha _2+6 \alpha _1 \alpha _2^2+3 \alpha _1 \alpha _4-6 \alpha _1^2 \alpha _4+\alpha
_4 \beta _1 &=0\\
1+5 \alpha _2+8 \alpha _2^2+6 \alpha _2^3+2 \alpha _4-2 \alpha _1 \alpha _4+3 \alpha _2 \alpha _4-6 \alpha _1 \alpha _2
\alpha _4 &=0\\
\alpha _4+2 \alpha _2 \alpha _4+2 \alpha _2^2 \alpha _4+\alpha _4^2-2 \alpha _1 \alpha _4^2&=0.
\end{array}
\end{equation}
Due to the equation 1 (of (3.2)) the following two cases occur:\\
\underline{Case 1. $\alpha _1+\alpha _1\alpha _2+\alpha _4\beta_1=0.$} Then
due to the equation 4 one has $\alpha _4=0$, in particular $-\alpha _1=\alpha
_1\alpha _2$, and $\alpha _1-2\alpha^2_2=0$. This implies that either
$\alpha_2=0$ or $\alpha_2=-1$. In both cases the equation 11 provides a
contradiction. So in Case 1 the system of equations (3.2) has no solution.\\
\underline{Case 2. $\beta_1=0.$} Then due to the equation 6 one has $\alpha_1=0, \ \frac{1}{2},\ \frac{1}{3}.$ In
$\alpha_1=0$ case the equation 7 provides a
contradiction, in $\alpha_1=\frac{1}{2}$ the equation 2 implies $\alpha_2=-1$
and the equation 4 implies $\alpha_4=0$, but the equation 3 provides a contradiction. In
$\alpha_1=\frac{1}{3}$ case
the equation 2 implies $\alpha_2=-\frac{2}{3}$ and the equation 4 implies $\alpha_4=0$,
but the equation 7 provides a contradiction. So there is no Jordan algebra among $A_1$.

For $A_{2}(\alpha_1, \beta_1, \beta_2)=\left(
\begin{array}{cccc}
\alpha _1 & 0 & 0 & 1 \\
\beta _1 & \beta _2 & -\alpha _1+1 & 0
\end{array}
\right)$  the system of equations (\ref{SE}) becomes

\[
\begin{array}{cccc}
-\beta _1^2 =0,\ \ \ \  -2 \beta _1+3 \alpha _1 \beta _1 =0,\ \ \ \  -1+2 \alpha _1-\alpha _1^2+\beta _2^2 =0,\ \ \ \  -\alpha _1 \beta
_1+\alpha _1^2 \beta _1+\beta _1 \beta _2^2 =0\\
-\alpha _1+3 \alpha _1^2-2 \alpha _1^3-\alpha _1 \beta _2+\alpha _1^2 \beta _2+\beta _2^2-2 \alpha _1 \beta _2^2+\beta
_2^3 =0,\ \ \ \ -\beta _1+\alpha _1 \beta _1+\beta _1 \beta _2 =0\\
-1+3 \alpha _1-2 \alpha _1^2 =0,\ \ \ \  \alpha _1 \beta _1-\beta _1 \beta _2 =0,\ \ \ \ 2 \alpha _1-3 \alpha _1^2-\beta _2+2 \alpha
_1 \beta _2-\beta _2^2 =0,\ \ \ \ 1-2 \alpha _1 =0.
\end{array}
\]
Hence, one gets \[\beta_1=0,\ \  \alpha_1=\frac{1}{2} \ \  \mbox{and} \ \ \beta_2=\pm \frac{1}{2}\] and we obtain the following Jordan
algebras\\
\[J_1=A_2\left(\frac{1}{2},0,\frac{1}{2}\right)=\left(
\begin{array}{cccc}
\frac{1}{2} & 0 & 0 & 1 \\
0 & \frac{1}{2} & \frac{1}{2} & 0
\end{array}
\right),\ \ \ J_2=A_2\left(\frac{1}{2},0,-\frac{1}{2}\right)=\left(
\begin{array}{cccc}
\frac{1}{2} & 0 & 0 & 1 \\
0 & -\frac{1}{2} & \frac{1}{2} & 0
\end{array}
\right).\]

For $A_{3}(\beta_1, \beta_2)=\left(
\begin{array}{cccc}
0 & 1 & 1 & 0 \\
\beta _1 & \beta _2 & 1 & -1
\end{array}
\right)$  the equation 11 of  (\ref{SE}) gives a contradiction and therefore there is no Jordan algebra among $A_3$.

For $A_{4}(\alpha_1, \beta_2)=\left(
\begin{array}{cccc}
\alpha _1 & 0 & 0 & 0 \\
0 & \beta _2 & -\alpha _1+1 & 0
\end{array}
\right)$  the system of equations (\ref{SE}) is equivalent to \\
\[  -\alpha _1+3 \alpha _1^2-2 \alpha _1^3-\alpha _1 \beta _2+\alpha _1^2 \beta _2+\beta _2^2-2 \alpha _1 \beta _2^2+\beta
_2^3=0.\]
So  \[\beta _2=-1+2 \alpha _1 \ \ \mbox{ or } \ \ \beta _2=-\sqrt{\alpha _1-\alpha _1^2} \ \  \mbox{ or } \ \   \beta _2=\sqrt{\alpha _1-\alpha
	_1^2}.\] Note that:
\begin{itemize}
	\item in $\alpha_1=0,1 $ cases $\sqrt{\alpha _1-\alpha _1^2}=-\sqrt{\alpha _1-\alpha _1^2}=0,$
	\item in
	$\alpha _1=\frac{1}{10} \left(5 \pm\sqrt{5}\right)$ cases  $-1+2 \alpha _1$ is equal to one of $\pm\sqrt{\alpha
		_1-\alpha _1^2}.$
\end{itemize}
Therefore in this case one has the following Jordan algebras:\\
\[J_3(\alpha_1)=A_4(\alpha_1,-1+2\alpha_1)=\left(
\begin{array}{cccc}
\alpha _1 & 0 & 0 & 0 \\
0 & -1+2 \alpha _1 & -\alpha _1+1 & 0
\end{array}
\right),\ \ \mbox{where} \ \ \alpha _1\in \mathbb{F} \ \  \mbox{and}\ \ \alpha _1\neq\frac{1}{10} \left(5\pm\sqrt{5}\right),\]
\[J_4=A_4\left(\alpha_1,\sqrt{\alpha _1-\alpha _1^2}\right)=\left(
\begin{array}{cccc}
\alpha _1 & 0 & 0 & 0 \\
0 & \sqrt{\alpha _1-\alpha _1^2} & -\alpha _1+1 & 0
\end{array}
\right),\ \ \mbox{where} \ \ \alpha _1\in \mathbb{F},\]  \[J_5=A_4\left(\alpha_1,-\sqrt{\alpha _1-\alpha _1^2}\right)=\left(
\begin{array}{cccc}
\alpha _1 & 0 & 0 & 0 \\
0 & -\sqrt{\alpha _1-\alpha _1^2} & -\alpha _1+1 & 0
\end{array}
\right),\ \ \mbox{where} \ \ \alpha _1\in \mathbb{F}\ \ \mbox{and} \ \ \alpha _1\neq 0,1.\]

For $A_{5}(\alpha_1)=\left(
\begin{array}{cccc}
\alpha _1 & 0 & 0 & 0 \\
1 & 2\alpha _1-1 & -\alpha _1+1 & 0
\end{array}
\right)$  the system of equations (\ref{SE}) is equivalent to\\
\[ 1-5 \alpha _1+5 \alpha _1^2=0\] and therefore one has\\
\[J_6=A_5\left(\frac{1}{10} \left(5-\sqrt{5}\right)\right)=\left(
\begin{array}{cccc}
\frac{1}{10} \left(5-\sqrt{5}\right) & 0 & 0 & 0 \\
1 & -\frac{\sqrt{5}}{5} & \frac{1}{10} \left(5+\sqrt{5}\right) & 0
\end{array}
\right),\] \[J_7=A_5\left(\frac{1}{10} \left(5+\sqrt{5}\right)\right)=\left(
\begin{array}{cccc}
\frac{1}{10} \left(5+\sqrt{5}\right) & 0 & 0 & 0 \\
1 & \frac{\sqrt{5}}{5} & \frac{1}{10} \left(5-\sqrt{5}\right) & 0
\end{array}
\right).\]

For $A_{6}(\alpha_1, \beta_1)=\left(
\begin{array}{cccc}
\alpha _1 & 0 & 0 & 1 \\
\beta _1 & -\alpha _1+1 & -\alpha _1 & 0
\end{array}
\right)$ we obtain \\
$\begin{array}{cccc}
-\beta _1^2 =0, & 3 \alpha _1 \beta _1 =0, & 1-2 \alpha _1 =0, & \beta _1-2 \alpha _1 \beta _1+2 \alpha _1^2 \beta _1 =0, \\
1-5 \alpha _1+8 \alpha _1^2-6 \alpha _1^3 =0, & -2 \alpha _1^2 =0, & -\beta _1+2 \alpha _1 \beta _1=0, & -1+4 \alpha _1-6
\alpha _1^2 =0
\end{array}$\\ and evidently it has no solution.

For $A_{7}(\beta_1)=\left(
\begin{array}{cccc}
0 & 1 & 1 & 0 \\
\beta _1 & 1 & 0 & -1
\end{array}
\right)$ case the equation 2 of (\ref{SE}) gives a contradiction.

For $A_{8}(\alpha _1)=\left(
\begin{array}{cccc}
\alpha _1 & 0 & 0 & 0 \\
0 & 1-\alpha _1 & -\alpha _1 & 0
\end{array}
\right)$  the system of equations (\ref{SE}) is equivalent to\\
\[ 1-5 \alpha _1+8 \alpha _1^2-6 \alpha _1^3=0,\] and therefore
\[\alpha _1=\frac{1}{3} \ \ \mbox{or} \ \ \alpha _1= \frac{1}{2}-\frac{i}{2}\ \  \mbox{or}\ \ \alpha _1=\frac{1}{2}+\frac{i}{2},\] where $i$ is a
fixed element such that $i^2=-1$. Due to it one gets Jordan algebras
\[J_8=A_8\left(\frac{1}{3}\right)=\left(
\begin{array}{cccc}
\frac{1}{3} & 0 & 0 & 0 \\
0 & \frac{2}{3} & -\frac{1}{3} & 0
\end{array}
\right),\ J_9=A_8\left(\frac{1}{2}-\frac{i}{2}\right)=\left(
\begin{array}{cccc}
\frac{1}{2}-\frac{i}{2} & 0 & 0 & 0 \\
0 & \frac{1}{2}+\frac{i}{2} & -\frac{1}{2}+\frac{i}{2} & 0
\end{array}
\right),\] \[J_{10}=A_8\left(\frac{1}{2}+\frac{i}{2}\right)=\left(
\begin{array}{cccc}
\frac{1}{2}+\frac{i}{2} & 0 & 0 & 0 \\
0 & \frac{1}{2}-\frac{i}{2} & -\frac{1}{2}-\frac{i}{2} & 0
\end{array}
\right).\]

In $A_9,\ A_{10},\ \ A_{11}$ cases there are no Jordan algebras.

The algebra $A_{12}=\left(
\begin{array}{cccc}
0 & 0 & 0 & 0 \\
1 & 0 & 0 & 0
\end{array}
\right)$ is a Jordan algebra.
%\end{proof}

Accordingly, for $Char(\mathbb{F})=2,3$ and $5$ cases we state the following theorems.
\begin{theorem}{\label{thm2}}
	Over an algebraically closed field $\mathbb{F}$ ($Char(\mathbb{F})=2$) any nontrivial $2$-dimensional Jordan algebra is
	isomorphic to only one of the following algebras listed by their matrices of structure constants:
	\begin{itemize}
		\item $J_{1,2}=A_{3,2}(0,0)=\left(
		\begin{array}{cccc}
		0 & 1 & 1 & 0 \\
		0 & 0 & 1 & 1
		\end{array}
		\right),$ 
		\item $J_{2,2}=A_{3,2}(1,0)=\left(
		\begin{array}{cccc}
		1 & 1 & 1 & 0 \\
		0 & 0 & 0 & 1
		\end{array}
		\right),$
		\item $J_{3,2}(\alpha_1)=A_{4,2}(\alpha_1,1)=\left(
		\begin{array}{cccc}
		\alpha _1 & 0 & 0 & 0 \\
		0 & 1 & \alpha _1+1 & 0
		\end{array}
		\right),$  where $\alpha_1\in\mathbb{F},$
		\item $J_{4,2}(\alpha_1)=A_{4,2}\left(\alpha_1,\sqrt{\alpha _1+\alpha _1^2}\right)=\left(
		\begin{array}{cccc}
		\alpha _1 & 0 & 0 & 0 \\
		0 & \sqrt{\alpha _1+\alpha _1^2} & \alpha _1+1 & 0
		\end{array}
		\right),$\\ \hspace*{\fill}  where $\alpha_1\in\mathbb{F}$ and $\alpha _1^2+\alpha_1+1\neq 0,$
		\item $J_{5,2}=A_{5,2}(\alpha_1)=\left(
		\begin{array}{cccc}
		\alpha _1 & 0 & 0 & 0 \\
		1 & 1 & \alpha _1+1 & 0
		\end{array}
		\right),$ where $\alpha_1\in\mathbb{F}$ and $\alpha _1^2+\alpha_1+1=0,$
		\item $J_{6,2}=A_{8,2}(1)=\left(
		\begin{array}{cccc}
		1 & 0 & 0 & 0 \\
		0 & 0 & 1 & 0
		\end{array}
		\right),$
		\item $J_{7,2}=A_{10,2}=\left(
		\begin{array}{cccc}
		0 & 1 & 1 & 0 \\
		0 & 0 & 0 & 1
		\end{array}
		\right),$
		\item $ J_{8,2}=A_{12,2}=\left(
		\begin{array}{cccc}
		0 & 0 & 0 & 0 \\
		1 & 0 & 0 & 0
		\end{array}
		\right).$
	\end{itemize}
\end{theorem}
\begin{theorem}{\label{thm3}}
	Over an algebraically closed field $\mathbb{F}$ ($Char(\mathbb{F})=3$) any nontrivial $2$-dimensional Jordan algebra is
	isomorphic to only one of the following algebras listed by their matrices of structure constants, where $i\in\mathbb{F}$
	stands for a fixed element for which $i^2=-1$ :
	\begin{itemize}
		\item $J_{1,3}=A_{2,3}(-1,0,1)=\left(
		\begin{array}{cccc}
		-1 & 0 & 0 & 1 \\
		0 & 1 & -1 & 0
		\end{array}
		\right),$
		\item $J_{2,3}=A_{2,3}(-1,0,-1)=\left(
		\begin{array}{cccc}
		-1 & 0 & 0 & 1 \\
		0 & -1 & -1 & 0
		\end{array}
		\right),$
		\item $J_{3,3}=A_{4,3}(\alpha_1,-1-\alpha_1)=\left(
		\begin{array}{cccc}
		\alpha _1 & 0 & 0 & 0 \\
		0 & -1- \alpha _1 & -\alpha _1+1 & 0
		\end{array}
		\right),$\\ \hspace*{\fill} where $\alpha _1\in \mathbb{F}$ and $\alpha _1\neq -1\pm i,$
		\item $ J_{4,3}=A_{4,3}\left(\alpha_1,\sqrt{\alpha
			_1-\alpha _1^2}\right)=\left(
		\begin{array}{cccc}
		\alpha _1 & 0 & 0 & 0 \\
		0 & \sqrt{\alpha _1-\alpha _1^2} & -\alpha _1+1 & 0
		\end{array}
		\right),$
		\item $J_{5,3}=A_{4,3}\left(\alpha_1,-\sqrt{\alpha _1-\alpha _1^2}\right)=\left(
		\begin{array}{cccc}
		\alpha _1 & 0 & 0 & 0 \\
		0 & -\sqrt{\alpha _1-\alpha _1^2} & -\alpha _1+1 & 0
		\end{array}
		\right),$\\ \hspace*{\fill} where $\alpha _1\in \mathbb{F}$ and $\alpha _1\neq 0,1,$
		\item $J_{6,3}=A_{5,3}( - 1+i)=\left(
		\begin{array}{cccc}
		- 1+i & 0 & 0 & 0 \\
		1 & - i & -1-i & 0
		\end{array}
		\right),$
		\item $J_{7,3}=A_{5,3}(- 1-i)=\left(
		\begin{array}{cccc}
		- 1-i & 0 & 0 & 0 \\
		1 & i & - 1+i & 0
		\end{array}
		\right),$
		\item $J_{8,3}=A_{8,3}(-1+i)=\left(
		\begin{array}{cccc}
		-1+i & 0 & 0 & 0 \\
		0 & -1-i & 1-i & 0
		\end{array}
		\right),$
		\item $J_{9,3}=A_{8,3}(-1-i)=\left(
		\begin{array}{cccc}
		-1-i & 0 & 0 & 0 \\
		0 & -1+i & 1+i & 0
		\end{array}
		\right),$
		\item $J_{10,3}=A_{10,3}=\left(
		\begin{array}{cccc}
		0 & 1 & 1 & 0 \\
		0 & 0 & 0 & -1
		\end{array}
		\right),$
		\item $J_{11,3}=A_{12,3}=\left(
		\begin{array}{cccc}
		0 & 0 & 0 & 0 \\
		1 & 0 & 0 & 0
		\end{array}
		\right).$
	\end{itemize}
\end{theorem}
\begin{theorem}{\label{thm4}}
	Over an algebraically closed field $\mathbb{F}$ ($Char(\mathbb{F})=5$) any nontrivial $2$-dimensional Jordan algebra is
	isomorphic to only one of the following algebras listed by their matrices of structure constants, where $i\in\mathbb{F}$
	stands for a fixed element for which $i^2=-1$ :
	\begin{itemize}
		\item $J_{1,5}=A_2\left(\frac{1}{2},0,\frac{1}{2}\right)=\left(
		\begin{array}{cccc}
		\frac{1}{2} & 0 & 0 & 1 \\
		0 & \frac{1}{2} & \frac{1}{2} & 0
		\end{array}
		\right),$
		\item $J_{2,5}=A_2\left(\frac{1}{2},0,-\frac{1}{2}\right)=\left(
		\begin{array}{cccc}
		\frac{1}{2} & 0 & 0 & 1 \\
		0 & -\frac{1}{2} & \frac{1}{2} & 0
		\end{array}
		\right),$
		\item $J_{3,5}=A_4(\alpha_1,-1+2\alpha_1)=\left(
		\begin{array}{cccc}
		\alpha _1 & 0 & 0 & 0 \\
		0 & -1+2 \alpha _1 & -\alpha _1+1 & 0
		\end{array}
		\right),$
		\item $J_{4,5}=A_4\left(\alpha_1,\sqrt{\alpha _1-\alpha _1^2}\right)=\left(
		\begin{array}{cccc}
		\alpha _1 & 0 & 0 & 0 \\
		0 & \sqrt{\alpha _1-\alpha _1^2} & -\alpha _1+1 & 0
		\end{array}
		\right),$
		\item $J_{5,5}=A_4\left(\alpha_1,-\sqrt{\alpha _1-\alpha _1^2}\right)=\left(
		\begin{array}{cccc}
		\alpha _1 & 0 & 0 & 0 \\
		0 & -\sqrt{\alpha _1-\alpha _1^2} & -\alpha _1+1 & 0
		\end{array}
		\right),\ \mbox{where }\alpha _1\neq 0,1,$
		\item $J_{6,5}=A_8\left(\frac{1}{3}\right)=\left(
		\begin{array}{cccc}
		\frac{1}{3} & 0 & 0 & 0 \\
		0 & \frac{2}{3} & -\frac{1}{3} & 0
		\end{array}
		\right),$
		\item $J_{7,5}=A_8\left(\frac{1}{2}-\frac{i}{2}\right)=\left(
		\begin{array}{cccc}
		\frac{1}{2}-\frac{i}{2} & 0 & 0 & 0 \\
		0 & \frac{1}{2}+\frac{i}{2} & -\frac{1}{2}+\frac{i}{2} & 0
		\end{array}
		\right),$
		\item $J_{8,5}=A_8\left(\frac{1}{2}+\frac{i}{2}\right)=\left(
		\begin{array}{cccc}
		\frac{1}{2}+\frac{i}{2} & 0 & 0 & 0 \\
		0 & \frac{1}{2}-\frac{i}{2} & -\frac{1}{2}-\frac{i}{2} & 0
		\end{array}
		\right),$
		\item $J_{9,5}=A_9=\left(
		\begin{array}{cccc}
		\frac{1}{3} & 0 & 0 & 0 \\
		1 & \frac{2}{3} & \frac{-1}{3} & 0
		\end{array}
		\right)=\left(
		\begin{array}{cccc}
		-\frac{1}{2} & 0 & 0 & 0 \\
		1 & -1 & \frac{1}{2} & 0
		\end{array}
		\right),$
		\item $J_{10,5}=A_{12}=\left(
		\begin{array}{cccc}
		0 & 0 & 0 & 0 \\
		1 & 0 & 0 & 0
		\end{array}
		\right).$
	\end{itemize}
\end{theorem}

A Jordan algebra
$\mathbb{A}$ is said to be commutative if $\mathbf{u}\cdot \mathbf{v}=\mathbf{v}\cdot \mathbf{u}$ whenever
$\mathbf{u},\mathbf{v}\in \mathbb{A}$. From the results presented above one can easily derive the following
classification results on commutative Jordan algebras.
\begin{corollary}\label{C1}
	Over an algebraically closed field $\mathbb{F}$ of characteristic not $2,\ 3$ every nontrivial $2$-dimensional
	commutative Jordan algebra is isomorphic to only one algebra listed below by MSC:\\
	\[J_{1c}=\left(
	\begin{array}{cccc}
	\frac{1}{2} &0& 0 & 1 \\
	0 & \frac{1}{2} & \frac{1}{2} & 0
	\end{array}
	\right),\ J_{2c}=\left(
	\begin{array}{cccc}
	\frac{2}{3} & 0 & 0 & 0 \\
	0 & \frac{1}{3} & \frac{1}{3} & 0
	\end{array}
	\right),\ J_{3c}=\left(
	\begin{array}{cccc}
	1 & 0 & 0 & 0 \\
	0 & 0 & 0 & 0
	\end{array}
	\right),\] \[ J_{4c}=\left(
	\begin{array}{cccc}
	\frac{1}{2} & 0 & 0 & 0 \\
	0 & \frac{1}{2} & \frac{1}{2} & 0
	\end{array}
	\right),\ J_{5c}=\left(
	\begin{array}{cccc}
	0 & 0 & 0 & 0 \\
	1 & 0 & 0 & 0
	\end{array}
	\right).\]
\end{corollary}
\begin{corollary}\label{C2}
	Over any algebraically closed field $\mathbb{F}$ of characteristic $2$ every nontrivial $2$-dimensional commutative
	Jordan algebra is isomorphic to only one algebra listed below by MSC:\\
	\[J_{1c,2}=\left(
\begin{array}{cccc}
1 & 1 & 1 & 0 \\
0 & 0 & 0 & 1
\end{array}
\right),\
	J_{2c,2}=\left(
	\begin{array}{cccc}
	1 & 0 & 0 & 0 \\
	0 & 0 & 0 & 0
	\end{array}
	\right),\ J_{3c,2}=\left(
	\begin{array}{cccc}
	0 & 0 & 0 & 0 \\
	0 & 1 & 1 & 0
	\end{array}
	\right),\] \[J_{4c,2}=\left(
	\begin{array}{cccc}
	0 & 1 & 1 & 0 \\
	0 & 0 & 0 & 1
	\end{array}
	\right),\ J_{5c,2}=\left(
	\begin{array}{cccc}
	0 & 0 & 0 & 0 \\
	1 & 0 & 0 & 0
	\end{array}
	\right).\]
\end{corollary}
\begin{corollary}\label{C3}
	Over any algebraically closed field $\mathbb{F}$ of characteristic $3$ every nontrivial $2$-dimensional commutative
	Jordan algebra is isomorphic to only one algebra listed below by MSC:\\
	\[J_{1c,3}=\left(
	\begin{array}{cccc}
	-1 &0& 0 & 1 \\
	0 & -1 & -1 & 0
	\end{array}
	\right),\ J_{2c,3}=\left(
	\begin{array}{cccc}
	1 & 0 & 0 & 0 \\
	0 & 0 & 0 & 0
	\end{array}
	\right),\ J_{3c,3}=\left(
	\begin{array}{cccc}
	-1 & 0 & 0 & 0 \\
	0 & -1 & -1 & 0
	\end{array}
	\right),\] \[
	J_{4c,3}=\left(
	\begin{array}{cccc}
	0 & 1 & 1 & 0 \\
	0 & 0 & 0 & -1
	\end{array}
	\right), \ J_{5c,3}=\left(
	\begin{array}{cccc}
	0 & 0 & 0 & 0 \\
	1 & 0 & 0 & 0
	\end{array}
	\right).\]
\end{corollary}

\section{\bf Classification of $2$-dimensional real Jordan algebras}
Due to \cite{Bekbaev} we have the following classification theorem.
\begin{theorem}\label{th1} Any non-trivial 2-dimensional real algebra is isomorphic to only one of the following listed, by their matrices of structure constants, algebras:\\
	$A_{1,r}(\mathbf{c})=\left(
	\begin{array}{cccc}
	\alpha_1 & \alpha_2 &\alpha_2+1 & \alpha_4 \\
	\beta_1 & -\alpha_1 & -\alpha_1+1 & -\alpha_2
	\end{array}\right),\ \mbox{where}\ \mathbf{c}=(\alpha_1, \alpha_2, \alpha_4, \beta_1)\in \mathbb{R}^4,$\\
	$A_{2,r}(\mathbf{c})=\left(
	\begin{array}{cccc}
	\alpha_1 & 0 & 0 & 1 \\
	\beta _1& \beta _2& 1-\alpha_1&0
	\end{array}\right), \mbox{where}\ \beta_1\geq 0,\ \mathbf{c}=(\alpha_1, \beta_1, \beta_2)\in \mathbb{R}^3,$\\
	$A_{3,r}(\mathbf{c})=\left(
	\begin{array}{cccc}
	\alpha_1 & 0 & 0 & -1 \\
	\beta _1& \beta _2& 1-\alpha_1&0
	\end{array}\right), \mbox{where}\ \beta_1\geq 0,\ \mathbf{c}=(\alpha_1, \beta_1, \beta_2)\in \mathbb{R}^3,$\\
	$A_{4,r}(\mathbf{c})=\left(
	\begin{array}{cccc}
	0 & 1 & 1 & 0 \\
	\beta _1& \beta _2 & 1&-1
	\end{array}\right),\ \mbox{where}\ \mathbf{c}=(\beta_1, \beta_2)\in \mathbb{R}^2,$\\
	$A_{5,r}(\mathbf{c})=\left(
	\begin{array}{cccc}
	\alpha _1 & 0 & 0 & 0 \\
	0 & \beta _2& 1-\alpha _1&0
	\end{array}\right),\ \mbox{where}\ \mathbf{c}=(\alpha_1, \beta_2)\in \mathbb{R}^2,$\\
	$A_{6,r}(\mathbf{c})=\left(
	\begin{array}{cccc}
	\alpha_1& 0 & 0 & 0 \\
	1 & 2\alpha_1-1 & 1-\alpha_1&0
	\end{array}\right),\ \mbox{where}\ \mathbf{c}=\alpha_1\in \mathbb{R},$\\
	$A_{7,r}(\mathbf{c})=\left(
	\begin{array}{cccc}
	\alpha_1 & 0 & 0 & 1 \\
	\beta _1& 1-\alpha_1 & -\alpha_1&0
	\end{array}\right), \mbox{where}\ \beta_1\geq 0,\ \mathbf{c}=(\alpha_1, \beta_1)\in \mathbb{R}^2,$\\
	$A_{8,r}(\mathbf{c})=\left(
	\begin{array}{cccc}
	\alpha_1 & 0 & 0 & -1 \\
	\beta _1& 1-\alpha_1 & -\alpha_1&0
	\end{array}\right), \mbox{where}\ \beta_1\geq 0,\ \mathbf{c}=(\alpha_1, \beta_1)\in \mathbb{R}^2,$\\
	$A_{9,r}(\mathbf{c})=\left(
	\begin{array}{cccc}
	0 & 1 & 1 & 0 \\
	\beta_1& 1& 0&-1
	\end{array}\right),\ \mbox{where}\ \mathbf{c}=\beta_1\in \mathbb{R},$\\
	$A_{10,r}(\mathbf{c})=\left(
	\begin{array}{cccc}
	\alpha_1 & 0 & 0 & 0 \\
	0 & 1-\alpha_1 & -\alpha_1&0
	\end{array}\right),\ \mbox{where}\ \mathbf{c}=\alpha_1\in\mathbb{R},$\\
	$ A_{11,r}=\left(
	\begin{array}{cccc}
	\frac{1}{3}& 0 & 0 & 0 \\
	1 & \frac{2}{3} & -\frac{1}{3}&0
	\end{array}\right),\
	\ A_{12,r}=\left(
	\begin{array}{cccc}
	0 & 1 & 1 & 0 \\
	1 &0&0 &-1
	\end{array}
	\right),$\\
	$ A_{13,r}=\left(
\begin{array}{cccc}
  0 & 1 & 1 & 0\\
  -1 &0&0 &-1
\end{array}
	\right),\
	\ A_{14,r}=\left(
\begin{array}{cccc}
  0 & 1 & 1 & 0 \\
  0 & 0 & 0 & -1
\end{array}
\right),\
	\ A_{15,r}=\left(
	\begin{array}{cccc}
0 & 0 & 0 & 0 \\
	1 &0& 0 &0
	\end{array}\right).$\end{theorem}
As an application of Theorem 4.1 we prove the next result.
\begin{theorem}{\label{theorem1}}
Any non-trivial two-dimensional real Jordan algebra is isomorphic to only one of the following listed, by their matrices of structure constants, algebras:\\
$J_{1,r}=A_{2,r}\left(\frac{1}{2},0,\frac{1}{2}\right)=\left(
\begin{array}{cccc}
 \frac{1}{2} & 0 & 0 & 1 \\
 0 & \frac{1}{2} & \frac{1}{2} & 0
\end{array}
\right),\\ J_{2,r}=A_{2,r}\left(\frac{1}{2},0,-\frac{1}{2}\right)=\left(
\begin{array}{cccc}
 \frac{1}{2} & 0 & 0 & 1 \\
 0 & -\frac{1}{2} & \frac{1}{2} & 0
\end{array}
\right),$\\ $J_{3,r}=A_{3,r}(\frac{1}{2},0,\frac{1}{2})=\left(
\begin{array}{cccc}
 \frac{1}{2} & 0 & 0 & -1 \\
 0 & \frac{1}{2} & \frac{1}{2} & 0
\end{array}
\right)$\\
$J_{4,r}=A_{3,r}(\frac{1}{2},0,-\frac{1}{2})=\left(
\begin{array}{cccc}
 \frac{1}{2} & 0 & 0 & -1 \\
 0 & -\frac{1}{2} & \frac{1}{2} & 0
\end{array}
\right)$\\$J_{5,r}(\alpha _1)=A_{5,r}(\alpha_1,-1+2\alpha_1)=\left(
\begin{array}{cccc}
 \alpha _1 & 0 & 0 & 0 \\
 0 & -1+2 \alpha _1 & -\alpha _1+1 & 0
\end{array}
\right),\ \mbox{where } \alpha _1\in \mathbb{R}$ and $\alpha _1\neq\frac{1}{10} \left(5\pm\sqrt{5}\right),$\\ $J_{6,r}(\alpha_1)=A_{5,r}\left(\alpha_1,\sqrt{\alpha _1-\alpha _1^2}\right)=\left(
\begin{array}{cccc}
 \alpha _1 & 0 & 0 & 0 \\
 0 & \sqrt{\alpha _1-\alpha _1^2} & -\alpha _1+1 & 0
\end{array}
\right),$\ where $0\leq \alpha _1\leq 1,$\\ $J_{7,r}(\alpha_1)=A_{5,r}\left(\alpha_1,-\sqrt{\alpha _1-\alpha _1^2}\right)=\left(
\begin{array}{cccc}
 \alpha _1 & 0 & 0 & 0 \\
 0 & -\sqrt{\alpha _1-\alpha _1^2} & -\alpha _1+1 & 0
\end{array}
\right),$ where $0< \alpha _1< 1,$\\ $J_{8,r}=A_{6,r}\left(\frac{1}{10}
\left(5-\sqrt{5}\right)\right)=\left(
\begin{array}{cccc}
 \frac{1}{10} \left(5-\sqrt{5}\right) & 0 & 0 & 0 \\
 1 & -\frac{\sqrt{5}}{5} & \frac{1}{10} \left(5+\sqrt{5}\right) & 0
\end{array}
\right),$\\ $J_{9,r}=A_{6,r}\left(\frac{1}{10} \left(5+\sqrt{5}\right)\right)=\left(
\begin{array}{cccc}
 \frac{1}{10} \left(5+\sqrt{5}\right) & 0 & 0 & 0 \\
 1 & \frac{\sqrt{5}}{5} & \frac{1}{10} \left(5-\sqrt{5}\right) & 0
\end{array}
\right),$\\ $J_{10,r}=A_{10,r}\left(\frac{1}{3}\right)=\left(
\begin{array}{cccc}
 \frac{1}{3} & 0 & 0 & 0 \\
 0 & \frac{2}{3} & -\frac{1}{3} & 0
\end{array}
\right),\ \ \ J_{11,r}=A_{15,r}=\left(
\begin{array}{cccc}
 0 & 0 & 0 & 0 \\
 1 & 0 & 0 & 0
\end{array}
\right).$
\end{theorem}
\textbf{Proof}. To prove we test the system of equations (\ref{SE}) for a consistency.

For $A_{1,r}(\alpha_1,\alpha_2,\alpha_4,\beta_1)=\left(
\begin{array}{cccc}
 \alpha _1 & \alpha _2 & \alpha _2+1 & \alpha _4 \\
 \beta _1 & -\alpha _1 & -\alpha _1+1 & -\alpha _2
\end{array}
\right)$ the system of equations (\ref{SE}) becomes (3.2).
 Due to the equation 1 (of (3.2)) the following two cases occur:\\
 \underline{Case 1. $\alpha _1+\alpha _1\alpha _2+\alpha _4\beta_1=0.$} Then
due to the equation 4 one has $\alpha _4=0$, in particular $-\alpha _1=\alpha
_1\alpha _2$, and $\alpha _1-2\alpha^2_2=0$. This implies that either
$\alpha_2=0$ or $\alpha_2=-1$. In both cases the equation 11 provides a
contradiction. So in Case 1 the system of equations (3.2) has no solution.\\
\underline{Case 2. $\beta_1=0.$} Due to the equation 6 one has $\alpha_1=0, \ \frac{1}{2},\ \frac{1}{3}.$ In
$\alpha_1=0$ case the equation 7 provides a
contradiction, in $\alpha_1=\frac{1}{2}$ the equation 2 implies $\alpha_2=-1$
and the equation 4 implies $\alpha_4=0$, but the equation 3 provides a contradiction. In
$\alpha_1=\frac{1}{3}$ case
the equation 2 implies $\alpha_2=-\frac{2}{3}$ and the equation 4 implies $\alpha_4=0$,
but the equation 7 provides a contradiction. So there is no Jordan algebra among $A_{1,r}$.

For $A_{2,r}(\alpha_1,\beta_1,\beta_2)=\left(
\begin{array}{cccc}
 \alpha _1 & 0 & 0 & 1 \\
 \beta _1 & \beta _2 & -\alpha _1+1 & 0
\end{array}
\right)$ the system of equations (\ref{SE}) becomes
\[   -\beta _1^2 =0,\   -2 \beta _1+3 \alpha _1 \beta _1 =0,\  -1+2 \alpha _1-\alpha _1^2+\beta _2^2 =0,\   -\alpha _1 \beta
 _1+\alpha _1^2 \beta _1+\beta _1 \beta _2^2 =0,\]
\[ -\alpha _1+3 \alpha _1^2-2 \alpha _1^3-\alpha _1 \beta _2+\alpha _1^2 \beta _2+\beta _2^2-2 \alpha _1 \beta _2^2+\beta
 _2^3 =0,\ -\beta _1+\alpha _1 \beta _1+\beta _1 \beta _2 =0,\]
 \[ -1+3 \alpha _1-2 \alpha _1^2 =0,\  \alpha _1 \beta _1-\beta _1 \beta _2 =0,\  2 \alpha _1-3 \alpha _1^2-\beta _2+2 \alpha
 _1 \beta _2-\beta _2^2 =0,\  1-2 \alpha _1 =0.
    \]
  Hence, one gets \[\beta_1=0,\ \  \alpha_1=\frac{1}{2} \ \  \mbox{and} \ \ \beta_2=\pm \frac{1}{2}\] and we obtain the following Jordan
  algebras\\
  $J_{1,r}=A_{2,r}\left(\frac{1}{2},0,\frac{1}{2}\right)=\left(
\begin{array}{cccc}
 \frac{1}{2} & 0 & 0 & 1 \\
 0 & \frac{1}{2} & \frac{1}{2} & 0
\end{array}
\right),\ \ \ J_{2,r}=A_{2,r}\left(\frac{1}{2},0,-\frac{1}{2}\right)=\left(
\begin{array}{cccc}
 \frac{1}{2} & 0 & 0 & 1 \\
 0 & -\frac{1}{2} & \frac{1}{2} & 0
\end{array}
\right)$.

For $A_{3,r}(\alpha_1,\beta_1,\beta_2)=\left(
	\begin{array}{cccc}
	\alpha_1 & 0 & 0 & -1 \\
	\beta _1& \beta _2& 1-\alpha_1&0
	\end{array}\right)$ the system of equations (\ref{SE}) becomes
\begin{equation}\label{SEb}\begin{array}{ccccc}
 \beta _1^2=0,\ \ \ \ 2 \beta _1-3 \alpha _1 \beta _1 =0,\ \ \ \  1-2 \alpha _1+\alpha _1^2-\beta _2^2=0,\ \ \ \ -\alpha _1 \beta _1+\alpha _1^2 \beta _1+\beta _1 \beta _2^2 =0,\\
 -\alpha _1+3 \alpha _1^2-2 \alpha _1^3-\alpha _1 \beta _2+\alpha _1^2 \beta _2+\beta _2^2-2 \alpha _1 \beta _2^2+\beta _2^3 =0,\ \ \ \ \ \beta _1-\alpha _1 \beta _1-\beta _1 \beta _2 =0,\\
 1-3 \alpha _1+2 \alpha _1^2=0, \ \ \ -\alpha _1 \beta _1+\beta _1 \beta _2 =0,\ \ \ -2 \alpha _1+3 \alpha _1^2+\beta _2-2 \alpha _1 \beta _2+\beta _2^2=0, \ \ \  1-2 \alpha _1=0.
\end{array}\end{equation}
Due to the equations $1$ and $3$ in the system of equations (\ref{SEb}) we get $\beta _1=0, \ \alpha _1=\frac{1}{2},$ then
$\beta _2=\pm \frac{1}{2},$ hence we obtain the following Jordan
  algebras\\
$J_{3,r}=A_{3,r}(\frac{1}{2},0,\frac{1}{2})=\left(
\begin{array}{cccc}
 \frac{1}{2} & 0 & 0 & -1 \\
 0 & \frac{1}{2} & \frac{1}{2} & 0
\end{array}
\right),\ J_{4,r}=A_{3,r}(\frac{1}{2},0,-\frac{1}{2})=\left(
\begin{array}{cccc}
 \frac{1}{2} & 0 & 0 & -1 \\
 0 & -\frac{1}{2} & \frac{1}{2} & 0
\end{array}
\right).$

For $A_{4,r}(\beta_1,\beta_2)=\left(
\begin{array}{cccc}
 0 & 1 & 1 & 0 \\
 \beta _1 & \beta _2 & 1 & -1
\end{array}
\right)$ the equation 11 of (\ref{SE}) gives a contradiction and therefore there is no Jordan algebra among $A_{4,r}$.

For $A_{5,r}(\alpha_1,\beta_2)=\left(
\begin{array}{cccc}
 \alpha _1 & 0 & 0 & 0 \\
 0 & \beta _2 & -\alpha _1+1 & 0
\end{array}
\right)$ the system of equations (\ref{SE}) is equivalent to \\
\[ -\alpha _1+3 \alpha _1^2-2 \alpha _1^3-\alpha _1 \beta _2+\alpha _1^2 \beta _2+\beta _2^2-2 \alpha _1 \beta _2^2+\beta
  _2^3=0.\]
  So $\beta _2$ equals one of $ -1+2 \alpha _1,\ -\sqrt{\alpha _1-\alpha _1^2},\
                \sqrt{\alpha _1-\alpha _1^2}$.
Note that in $\alpha_1=0,1 $ cases $\sqrt{\alpha _1-\alpha _1^2}=-\sqrt{\alpha _1-\alpha _1^2}=0$ and in
  $\alpha _1=\frac{1}{10} \left(5 \pm\sqrt{5}\right)$ cases $-1+2 \alpha _1$ is equal to one of $\pm\sqrt{\alpha
  _1-\alpha _1^2}$. Therefore in this case one has the following Jordan algebras\\
  $J_{5,r}(\alpha_1)=A_{5,r}(\alpha_1,-1+2\alpha_1)=\left(
\begin{array}{cccc}
 \alpha _1 & 0 & 0 & 0 \\
 0 & -1+2 \alpha _1 & -\alpha _1+1 & 0
\end{array}
\right),$\ 
where $\alpha _1\in \mathbb{R}$ and $ \alpha _1\neq\frac{1}{10} \left(5\pm\sqrt{5}\right),$\\
$J_{6,r}(\alpha_1)=A_{5,r}\left(\alpha_1,\sqrt{\alpha _1-\alpha _1^2}\right)=\left(
\begin{array}{cccc}
 \alpha _1 & 0 & 0 & 0 \\
 0 & \sqrt{\alpha _1-\alpha _1^2} & -\alpha _1+1 & 0
\end{array}
\right),$\  where $0\leq \alpha _1\leq 1,$\\ $J_{7,r}(\alpha_1)=A_{5,r}\left(\alpha_1,-\sqrt{\alpha _1-\alpha _1^2}\right)=\left(
\begin{array}{cccc}
 \alpha _1 & 0 & 0 & 0 \\
 0 & -\sqrt{ \alpha _1-\alpha _1^2} & -\alpha _1+1 & 0
\end{array}
\right),$\ where $0< \alpha _1 < 1.$

For $A_{6,r}(\alpha_1)=\left(
\begin{array}{cccc}
 \alpha _1 & 0 & 0 & 0 \\
 1 & 2\alpha _1-1 & -\alpha _1+1 & 0
\end{array}
\right)$ the system of equations (\ref{SE}) is equivalent to\\
\[ 1-5 \alpha _1+5 \alpha _1^2=0\] and therefore one has\\
    $J_{8,r}=A_{6,r}\left(\frac{1}{10} \left(5-\sqrt{5}\right)\right)=\left(
\begin{array}{cccc}
 \frac{1}{10} \left(5-\sqrt{5}\right) & 0 & 0 & 0 \\
 1 & -\frac{\sqrt{5}}{5} & \frac{1}{10} \left(5+\sqrt{5}\right) & 0
\end{array}\right),$\\ $J_{9,r}=A_{6,r}\left(\frac{1}{10} \left(5+\sqrt{5}\right)\right)=\left(
\begin{array}{cccc}
 \frac{1}{10} \left(5+\sqrt{5}\right) & 0 & 0 & 0 \\
 1 & \frac{\sqrt{5}}{5} & \frac{1}{10} \left(5-\sqrt{5}\right) & 0
\end{array}
\right).$

For $A_{7,r}(\alpha_1,\beta_1)=\left(
\begin{array}{cccc}
 \alpha _1 & 0 & 0 & 1 \\
 \beta _1 & -\alpha _1+1 & -\alpha _1 & 0
\end{array}
\right)$ we obtain \\
\[\beta _1^2 =0,\ \ \ \ \ 3 \alpha _1 \beta _1 =0,\ \ \ \ \ 1-2 \alpha _1 =0,\ \ \ \ \ \beta _1-2 \alpha _1 \beta _1+2 \alpha _1^2 \beta _1 =0,\ \ \ \ \
 1-5 \alpha _1+8 \alpha _1^2-6 \alpha _1^3 =0,\] \[-2 \alpha _1^2 =0,\ \ \ \ \  -\beta _1+2 \alpha _1 \beta _1=0,\ \ \ \ -1+4 \alpha _1-6
 \alpha _1^2 =0.\] and evidently it has no solution.

For $A_{8,r}(\alpha_1,\beta_1))=\left(
	\begin{array}{cccc}
	\alpha_1 & 0 & 0 & -1 \\
	\beta _1& 1-\alpha_1 & -\alpha_1&0
	\end{array}\right)$
case the equations 3 and 8 of (\ref{SE}) give a contradiction.

For $A_{9,r}(\beta_1)=\left(
\begin{array}{cccc}
 0 & 1 & 1 & 0 \\
 \beta _1 & 1 & 0 & -1
\end{array}
\right)$ case the equation 2 of (\ref{SE}) gives a contradiction.

For $A_{10,r}(\alpha_1)=\left(
\begin{array}{cccc}
 \alpha _1 & 0 & 0 & 0 \\
 0 & 1-\alpha _1 & -\alpha _1 & 0
\end{array}
\right)$ the system of equations (\ref{SE}) is equivalent to\\
\[ 1-5 \alpha _1+8 \alpha _1^2-6 \alpha _1^3=0,\] and therefore
   $\alpha _1=\frac{1}{3},$ and one gets Jordan algebra
  $J_{10,r}=A_{10,r}\left(\frac{1}{3}\right)=\left(
\begin{array}{cccc}
 \frac{1}{3} & 0 & 0 & 0 \\
 0 & \frac{2}{3} & -\frac{1}{3} & 0
\end{array}
\right).$

In $A_{11,r},\ A_{12,r},\ A_{13,r},\ \ A_{14,r}$ cases there are no Jordan algebras.
The algebra $A_{15,r}=\left(
\begin{array}{cccc}
 0 & 0 & 0 & 0 \\
 1 & 0 & 0 & 0
\end{array}
\right)$ is a Jordan algebra.
%\end{proof}
\begin{remark}
	In \cite{Wallace} one can find a classification of two-dimensional power associative real algebras given as $A_1, A_3, A(\sigma), B_1, B_2, B_3, B_4, B_5$, where $\sigma\in \mathbb{R}$. In our notations they are: \\
	$A_1$ is the trivial algebra, $A_3=\left(
	\begin{array}{cccc}
	1 & 0 & 0 & 0 \\
	0 & 0 & 1 & 0
	\end{array}
	\right)\simeq J_{5,r}(\frac{1}{2}),$\\ $A(\sigma)=\left(
	\begin{array}{cccc}
	1+\sigma & 0 & 0 & 0 \\
	0 & 1 & \sigma & 0
	\end{array}
	\right)\simeq J_{5,r}(\frac{1+\sigma}{1+2\sigma}),$ where $\sigma\neq -\frac{1}{2}, \frac{-1\pm\sqrt{5}}{2},$\\
	$A(\frac{-1}{2})=J_{10,r},$
	$A(\frac{-1+\sqrt{5}}{2})=J_{6,r}(\frac{5+\sqrt{5}}{10}),$\\
	$A(\frac{-1-\sqrt{5}}{2})=J_{7,r}(\frac{5-\sqrt{5}}{10}),$
	$ B_1=\left(
	\begin{array}{cccc}
	0 & 0 & 0 & 0 \\
	1 & 0 & 0 & 0
	\end{array}
	\right)\simeq J_{11,r},$\\ $B_2=\left(
	\begin{array}{cccc}
	0 & 0 & 0 & 0 \\
	1 & 1 & 1 & 1
	\end{array}
	\right)\simeq J_{6,r}(1),$  $B_3=\left(
	\begin{array}{cccc}
	0 & 1 & 1 & 0 \\
	1 & 0 & 0 & 1
	\end{array}
	\right)\simeq J_{1,r},$\\ $B_4=\left(
	\begin{array}{cccc}
	0 & 1 & 1 & 0 \\
	-1 & 0 & 0 & 1
	\end{array}
	\right)\simeq J_{3,r},$  $B_5=\left(
	\begin{array}{cccc}
	0 & 1 & 1 & 0 \\
	0 & 0 & 0 & 1
	\end{array}
	\right)\simeq J_{6,r}(\frac{1}{2})$.
	
	So the representatives of non power associative Jordan algebras are:\\
	$J_{2,r},\ J_{4,r},\ J_{6,r}(\alpha_1),$ where $0\leq \alpha_1 <1$ and $\alpha_1\neq \frac{1}{2}, \frac{5+\sqrt{5}}{10}$,\\ $J_{7,r}(\alpha_1),$ where $0< \alpha_1 <1$ and $\alpha_1\neq  \frac{5-\sqrt{5}}{10}$,\ $J_{8,r}$,\ $J_{9,r}.$\\ In particular,
	we would like to note that $J_{6,r}(0)=\left(
	\begin{array}{cccc}
		0 & 0 & 0 & 0 \\
		0 & 0 & 1 & 0
	\end{array}
	\right)$ is an example of Jordan algebra which is not power associative over any field $\mathbb{F}$.\end{remark}

 From the results presented above one can easily derive the following
classification results on commutative Jordan algebras.
\begin{corollary}\label{T1}
  Any non-trivial two-dimensional commutative real Jordan algebra is isomorphic to only one of the following listed, by their matrices of structure constants, algebras:\\
 $J_{1c,r}=\left(
\begin{array}{cccc}
 \frac{1}{2} &0& 0 & 1 \\
 0 & \frac{1}{2} & \frac{1}{2} & 0
\end{array}
\right),\ J_{2c,r}=\left(
\begin{array}{cccc}
 \frac{1}{2} &0& 0 & -1 \\
 0 & \frac{1}{2} & \frac{1}{2} & 0
\end{array}
\right),\ J_{3c,r}=\left(
\begin{array}{cccc}
 \frac{2}{3} & 0 & 0 & 0 \\
 0 & \frac{1}{3} & \frac{1}{3} & 0
\end{array}
\right),\\ J_{4c,r}=\left(
\begin{array}{cccc}
 1 & 0 & 0 & 0 \\
 0 & 0 & 0 & 0
\end{array}
\right),$\ $J_{5c,r}=\left(
\begin{array}{cccc}
 \frac{1}{2} & 0 & 0 & 0 \\
 0 & \frac{1}{2} & \frac{1}{2} & 0
\end{array}
\right),\ J_{6c,r}=\left(
\begin{array}{cccc}
 0 & 0 & 0 & 0 \\
 1 & 0 & 0 & 0
\end{array}
\right).$
\end{corollary}
\begin{remark} Classification of nontrivial two-dimensional real commutative Jordan algebras has been presented in \cite{Wallace} as the algebras $A(1), B_1, B_2, B_3, B_4, B_5$ and they are isomorphic respectively to the algebras $J_{3c,r}, J_{6c,r}, J_{4c,r}, J_{1c,r}, J_{2c,r}, J_{5c,r}$ from Corollary 3.3.\end{remark}
	
 \begin{remark} When we compare our results with that of \cite{Bermúdez} we find that the system of equations (\ref{SE}) in commutative case ($\alpha_3=\alpha_2,\ \beta_3=\beta_2$) and the system of equations (4) in \cite{Bermúdez} are the same. There are the following isomorphisms between representatives of Corollary 3.3 and that of \cite{Bermúdez} :\\
  $J_{1c,r} \simeq \psi_0,$ $J_{2c,r}\simeq \psi_5,$ $J_{3c,r} \simeq \psi_4,$ $J_{4c,r} \simeq \psi_2,$ $J_{5c,r} \simeq \psi_1,$ and $J_{6c,r} \simeq \psi_3.$
\end{remark}

\end{document}